\documentclass[10pt,notitlepage,twoside,a4paper]{amsart}
 \usepackage{amsfonts}

\usepackage{amsmath,amssymb,enumerate}

\usepackage{epsfig,fancyhdr,color}

\usepackage{amssymb}
\usepackage{amsmath,amsthm}
\usepackage{latexsym}
\usepackage{amscd}
\usepackage{psfrag}
\usepackage{graphicx}
\usepackage[latin1]{inputenc}
\usepackage[all]{xy}
\usepackage[mathcal]{eucal}

\renewcommand{\textsc}{\textcolor{red}}

\newtheorem{theorem}{\rm\bf Theorem}[section]
\newtheorem{proposition}[theorem]{\rm\bf Proposition}
\newtheorem{lemma}[theorem]{\rm\bf Lemma}

\newtheorem*{theorem 1}{\rm\bf Proposition 1}
\newtheorem*{theorem 2}{\rm\bf Proposition 2}

\theoremstyle{definition}
\newtheorem{definition}[theorem]{\rm\bf Definition}

\theoremstyle{remark}
\newtheorem{remark}[theorem]{\rm\bf Remark}

\newtheorem{question}[theorem]{\rm\bf Question}

\def\interieur#1{\mathord{\mathop{\kern 0pt #1}\limits^\circ}}

\title[Length spectra and the Teichm\"uller metric]{Length spectra and the Teichm\"uller metric for surfaces with boundary}

\author{Lixin Liu}
\address{Lixin Liu, Department of Mathematics, Zhongshan University, 510275, Guangzhou, P. R. China}
\email{mcsllx@mail.sysu.edu.cn}

\author{Athanase Papadopoulos}
\address{Athanase Papadopoulos, Max-Plank-Institut f\"ur Mathematik, Vivatsgasse 7, 53111 Bonn, Germany and: Institut de Recherche Math{\'e}matique Avanc\'ee,
Universit{\'e} de Strasbourg and CNRS,
7 rue Ren\'e Descartes,
 67084 Strasbourg Cedex, France} \email{papadopoulos@math.u-strasbg.fr}
\date{\today}

\author{Weixu Su}
\address{Weixu Su, Department of Mathematics, Zhongshan University, 510275, Guangzhou, P. R. China}
\email{su023411040@163.com}

\author{Guillaume Th\'eret}
\address{Guillaume Th\'eret, Max-Plank-Institut f\"ur Mathematik, Vivatsgasse 7, 53111 Bonn, Germany}
\email{theret@mpim-bonn.mpg.de}

\begin{document}

\begin{abstract} 

We consider some metrics and weak metrics defined on the Teichm\"uller space of 
 a surface of finite type with nonempty boundary, that are defined 
using the hyperbolic length spectrum of simple closed curves and of properly embedded arcs, and we compare these metrics and weak metrics with the Teichm\"uller metric. The comparison is on subsets of Teichm\"uller space which we call ``$\varepsilon_0$-relative $\epsilon$-thick parts", 
and whose definition depends on the choice of some positive constants $\varepsilon_0$ and $\epsilon$.
Meanwhile, we give a formula for the Teichm\"uller metric of a surface with boundary in terms of extremal lengths of families of arcs.

\bigskip

\medskip

\medskip

\noindent Lixin Liu and Weixu Su were partially supported by NSFC (No. 10871211).
 
\end{abstract}
\maketitle

\tableofcontents

\section{Introduction}\label{intro}

In this paper, $S$ is a connected oriented surface of finite topological type with nonempty boundary. 
The surface $S$ is obtained from a closed oriented surface by removing a finite number $p\geq 0$ of punctures 
and a finite number $b\geq 1$ of disjoint open disks. 
We shall say that $S$ has $p$ punctures and $b$ boundary components, or that it is of type $(g,p,b)$.
The boundary of $S$ is denoted by $\partial S$, and we shall assume throughout the paper that $\partial S\not=\emptyset$ 
unless the contrary is explicitly stated.
The Euler characteristic of $S$ is equal to $\chi(S)= 2-2g-p-b$, and we shall always assume that $\chi(S)<0$.\\

We shall consider conformal (or Riemann surface) structures on $S$. By the expression ``conformal structure"  on a surface $S$ with nonempty boundary, we mean that the topological interior of $S$ is equipped with a conformal structure in the usual sense, and
 that the following two properties are satisfied:
\begin{enumerate}
\item each puncture has a neighborhood which is conformally equivalent to the punctured disk $\{z\in \mathbb{C}\  | \ 0 < |z|< 1 \}$;
\item each boundary component has a neighborhood which is conformally equivalent to an
annulus of finite modulus $\{z\in \mathbb{C} \  | \  a \leq  |z| < 1\}$ for some positive constant $a$ that depends on the puncture. 
\end{enumerate}

To simplify notation, we shall also denote by the same letter $S$ the surface equipped with a conformal structure. 
Then, $\bar{S}$ denotes the surface $S$ equipped with the {\it mirror-image} conformal structure. 
We recall that this means that a conformal atlas for $\bar{S}$ is obtained from an atlas for $S$ by composing each local coordinate $z$ for $S$ 
with the conjugation map $z\mapsto \bar{z}$. 
We let $S^d= S\cup \bar{S}$ denote the {\it conformal double of $S$}, 
obtained by gluing $S$ and $\bar{S}$ along corresponding boundary components, using the identity map. 
The surface $S^d$ is equipped with a natural anti-holomorphic involution, whose fixed point set is the boundary of $S$, considered as embedded in $S^d$. 
By the Poincar\'e uniformization theorem, $S^d$ carries a unique complete hyperbolic metric which is compatible with its conformal structure. 
Because of the symmetry in the definition of $\bar{S}$, this hyperbolic metric is invariant under the natural anti-holomorphic involution of $S^d$. 
The {\it canonical hyperbolic metric} on $S$ is the restriction to $S$ (with respect to the natural inclusion $S\subset S^d$) of this complete hyperbolic metric on $S^d$.  
The conformal structure that underlies the canonical hyperbolic metric on $S$ is the conformal structure that we started with.

We denote by $\mathcal {T}(S)$ the reduced Teichm\"{u}ller space of marked conformal structures on $S$. 
This is the set of equivalence classes of pairs $(X, f)$, where $X$ is a Riemann surface and $f : S\rightarrow X$ is a
homeomorphism (called the marking), and where $(X_1, f_1)$ is considered to be equivalent to $(X_2, f_2)$ if there is a conformal map 
$h : X_1\rightarrow X_2$ which is homotopic to $f_2\circ f_1^{-1}$. 
We recall that in the reduced theory, homotopies need not fix the boundary of $S$ pointwise. 
Since all Teichm\"uller spaces that we consider are reduced, we shall omit the word ``reduced" in our exposition. 
Furthermore, we shall often denote an element $(X, f)$ of
$\mathcal {T}(S)$ by $X$, without explicit reference to the marking.

The elements of $\mathcal {T}(S)$ can also be taken to be (equivalence classes of) hyperbolic metrics with geodesic boundary, 
and these metrics will always be the canonical hyperbolic metrics for the conformal structures they induce.

In the paper \cite{LPST1}, we studied some metrics and weak metrics on the Teichm\"uller space $\mathcal {T}(S)$, 
defined using  the hyperbolic length spectrum of simple closed curves and of properly embedded arcs in $S$. 
In the present paper, we compare these metrics and weak metrics with the Teichm\"uller metric on $\mathcal {T}(S)$, 
on subsets of this space which we call ``$\varepsilon_0$-relative $\epsilon$-thick parts", 
whose definition depends on the choice of some positive constants $\varepsilon_0$ and $\epsilon$, and which we describe below.

\section{On the geometry of surfaces with boundary}\label{section:conformal}

We need to give some complements to Thurston's theory concerning a surface $S$ with boundary. 
As a general rule, the point is to include in the theory the study of arcs joining boundary components, instead of dealing only with simple closed curves on the surface.  
As we shall see, many results on spaces associated to a surface $S$ with boundary can be recovered from corresponding results on a surface without boundary 
by taking the double of $S$. 
In this section, we recall some definitions and give precise statements concerning some topological and conformal notions on surfaces with boundary.
Abikoff's book \cite{Abikoff} deals with surfaces with boundary, 
and there is also a small subsection on quadratic differentials on surfaces with boundary in Strebel's book (see \cite{Strebel} p. 157).\\

We first recall a few topological definitions.
A simple closed curve in $S$ is said to be {\it peripheral} if it is homotopic to a puncture. 
It is said to be {\it essential} if it is non-peripheral and if it is not homotopic to a point (but it can be homotopic to a boundary component).

We let $\mathcal{C}=\mathcal{C}(S)$ be the set of homotopy classes of essential simple closed curves on $S$. 

An {\it arc} in $S$ is the homeomorphic image of a closed interval which is properly embedded in $S$ 
(that is, the interior of the arc is in the interior of $S$ and the endpoints of the arc are on the boundary of $S$). 
All homotopies of arcs that we consider are relative to $\partial S$, that is, they leave the endpoints of arcs on the set $\partial S$ 
(but they do not necessarily fix pointwise the endpoints).
An arc is said to be {\it essential} if it is not homotopic (relative to $\partial S$) to a map whose image is in $\partial S$. 

We let $\mathcal{B}=\mathcal{B}(S)$ be the union of the set of homotopy classes
of essential arcs on $S$ with the set of homotopy classes of simple closed curves which are homotopic to boundary components. 
 
We shall denote by $[\gamma]\in \mathcal{B}\cup\mathcal{C}$ the equivalence class of a simple closed curve or of an arc $\gamma$ on $S$.

For any $[\gamma]\in\mathcal{B}\cup\mathcal{C}$ and for any $X\in \mathcal{T}(S)$, 
we let $\gamma^X$ be the geodesic representative of $[\gamma]$ 
(that is, the curve of shortest length in the homotopy class relative to $\partial S$) 
with respect to the hyperbolic structure $X$, if $X$ is a hyperbolic structure, or with respect to the canonical hyperbolic metric
associated to the conformal structure $X$, if $X$ is only a conformal structure. 
The geodesic $\gamma^X$ is unique; it coincides with a boundary component if $[\gamma]$ is the homotopy class of that boundary component, and it is
orthogonal to $\partial S$ at each intersection point if $[\gamma]$ is the homotopy class of an arc. 
We denote by $l_X(\gamma)$ the length of $\gamma^X$ with respect to the hyperbolic metric $X$.
To simplify notation, we shall sometimes denote $[\gamma]$ by $\gamma$.

We denote by $\mathcal{ML}(S)$ the space of measured geodesic laminations on $S$. 
This space is equipped with a topology defined by Thurston (cf. \cite{Thurston-notes}).

We recall that any measured geodesic lamination on a hyperbolic surface of finite topological type can be decomposed as the union of finitely many components 
whose support is of one of the following three types:
 \begin{enumerate}
 \item an essential geodesic arc;
 \item a boundary component;
 \item a geodesic lamination in the interior of $S$ which is minimal, that is, in which every leaf is dense.
 \end{enumerate}
 
Finally, we recall that there are natural homeomorphisms between the spaces $\mathcal{ML}(S)$, when the hyperbolic structure on $S$ varies, 
so that it is possible to talk about a measured geodesic laminations on $S$ without referring to a specific hyperbolic structure on that surface. 
Such a space  $\mathcal{ML}(S)$ is also canonically homeomorphic to the space $\mathcal{MF}(S)$ of equivalence classes of measured foliations on $S$.
An arc (respectively, hyperbolic structure, simple closed curve, measured lamination, etc.) on $S^d$ is said to be {\it symmetric} 
if it is invariant by the natural involution.\\

\begin{figure}[!hbp]
\centering
\includegraphics[width=.4\linewidth]{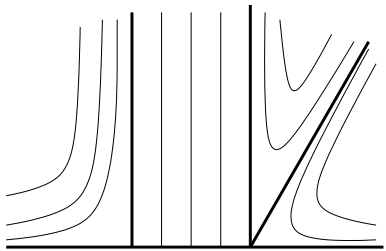}
\caption{\small{The leaves of a measured foliation can be either transverse or tangent to the boundary component of the surface.}}
\label{boundaryfoliation}
\end{figure}

In what follows, quadratic differentials are as usual assumed to be meromorphic with at most simple poles at punctures.
Quadratic differentials on a conformal surface with boundary are assumed to take real values on vectors tangent to the boundary components. 
More precisely, such a quadratic differential $\phi$ can have (isolated) zeroes on $\partial S$, 
and it satisfies $\phi(z)dz^2>0$ or $\phi(z)dz^2<0$ for any tangent vector $dz$ which is  tangent to $\partial S$  
at a point on  $\partial S$ which is in the complement of the zeroes of $\phi$.

A quadratic differential equips the underlying surface with a pair of measured foliations, namely, its horizontal and vertical foliations. 
On surfaces with boundary, measured foliations are allowed to be transverse to the boundary or tangent to the boundary, 
and both cases can occur at the same boundary component, as in Figure  \ref{boundaryfoliation}.  
By a result of Hubbard and Masur (cf. \cite{HM}, adapted to the case of surfaces with boundary), 
the ordered pair consisting of the horizontal and vertical measured foliation of a quadratic differential  
(and even, the pair of their Whitehead equivalence classes in the sense of Thurston) completely determines this quadratic differential.

Given a surface $S$ equipped with a quadratic differential $\phi$, 
each point on a boundary component of $S$ which is a nonzero point of $\phi$ has a neighborhood where the horizontal and vertical measured foliations 
are either perpendicular or tangent to the boundary. 
(As already said, both kinds of behavior can occur at a given boundary component, as in Figure \ref{boundaryfoliation}.) 
The quadratic differential $\phi$ gives, by reflection along $\partial S$, a well-defined quadratic differential on the double $S^d$ of $S$.

The Teichm\"{u}ller metric $d_{T}$ on $\mathcal{T}(S)$ is defined by
$$d_{T}(X,Y)=\frac{1}{2} \inf  \log K(f)$$
for any two marked conformal structures $X$ and $Y$ on $S$,  
where the infimum is taken over all quasiconformal homeomorphisms $f:X\to Y$ and where 
$K(f)$ is the dilatation of $f$.

If $X$ and $Y$ are two conformal structures on $S$, 
then a quasiconformal homeomorphism $f:X\to Y$ extends to a quasiconformal homeomorphism $f^d:X^d\to Y^d$ between the conformal doubles, 
with the same dilatation $K(f)$.

A Teichm\"uller map $f:X\to Y$ extends to a Teichm\"uller map $f^d:X^d\to Y^d$. 
From the uniqueness of Teichm\"uller maps, $f^d$ commutes with the natural anti-holomorphic involutions of the Riemann surfaces $X^d$ and $Y^d$, 
and its  associated initial quadratic differential is the double of the quadratic differential of $f$. 
In particular, the horizontal and vertical measured foliations of $f^d$ are symmetric with respect to the natural involution of $X^d$.

\begin{theorem}[Wolpert's inequality for surfaces with boundary]\label{th:Wolpert}
Let $X$ and $Y$ be two hyperbolic structures on $S$ and let $f:X\to Y$ be a
$K$-quasiconformal homeomorphism. For any element $\gamma$ which is either an essential
simple closed curve or an essential arc on $S$, we have
 \[\frac{l_Y(f(\gamma))}{l_X(\gamma)}\leq K.\]
\end{theorem}

The result for closed curves on surfaces without boundary is in Wolpert \cite{Wolpert}. 
The result for arcs follows by doubling. 
Indeed, doubling a quasiconformal homeomorphism produces on the doubled surface a quasiconformal homeomorphism with the same dilatation. 
Doubling a geodesic $X$- or $Y$-arc produces a simple closed geodesic on the doubled hyperbolic surface $X^d$ or $Y^d$ respectively. 
Wolpert's inequality for arcs in $S$ follows then from the inequality for simple closed curves on the double.

\section{Extremal length and the Teichm\"uller metric}

The notion of extremal length for families of curves or of arcs on a Riemann surface was introduced by Ahlfors and Beurling \cite{Beurling-Ahlfors}.
We follow the exposition that Kerckhoff gives in \cite{K} for elements of $\mathcal{C}(S)$, which we adapt to elements of $\mathcal{B}(S)\cup\mathcal{C}(S)$. 

There are two equivalent definitions of extremal length, one which Kerckhoff calls analytic, and another one which he calls geometric.

We consider a fixed Riemann surface structure on $S$. 
A {\it conformal metric} on $S$ is a Riemannian metric which may have isolated singularities 
and whose length element in every holomorphic coordinate $z$ in the complement of the singularities can be written as $\rho(z)\vert dz\vert$, 
where $\rho$ is a nonnegative real-valued function.  
We let $A_X$ denote the total area of the metric $X$. 

\begin{definition}[Analytic] 
The extremal length of $\gamma\in \mathcal{C}(S) \cup \mathcal{B}(S)$ is the quantity
\[E_{X}(\gamma)=\sup_{X}\frac{l^{2}_X(\gamma)}{A_X}\]
where the supremum is taken over all conformal metrics satisfying $0<A_X<\infty$.
\end{definition}

Given a conformal metric $X$ on $S$, we shall also use the notation $E_X(\gamma)$ to denote the extremal length of a family $\gamma$ 
of closed curves or of arcs on $S$, meaning that we consider the extremal length of $\gamma$ with respect to the conformal structure underlying $X$. 
In general, the metrics $X$ that we shall consider will be hyperbolic.

Since the area of any hyperbolic metric $X$ on $S$ is equal to $2\pi |\chi(S)|$, 
it follows from the  above definition that for any element $\gamma$ in $\mathcal{B}\cup\mathcal{C}$ we have
\begin{equation}\label{eq:extremal}
E_X(\gamma) \geq \frac{l_X^2(\gamma)}{2\pi |\chi(S)|}.
\end{equation}

For the geometric definition of conformal length, we need to consider cylinders and quadrilaterals in Riemann surfaces, and we recall the definitions.

A {\it cylinder} in a Riemann surface $S$ is a surface with boundary which is homeomorphic to a compact Euclidean annulus, which is immersed in $S$ 
in such a way that its interior is embedded. 
If $C$ is any closed curve in $S$ homotopic to a boundary component of a cylinder, then we say that $C$ is a {\it core curve} of the cylinder.

We shall consider quadrilaterals only on Riemann surfaces with nonempty boundary. 
A {\it quadrilateral} in such a surface $S$ is a closed disk with two disjoint distinguished curves on its boundary, the disk being immersed in 
$S$ with embedded interior, such that the distinguished curves are arcs in $S$ in the sense defined above, 
and such that the two curves in the complement of the distinguished curves in the boundary of the quadrilateral are contained in $\partial S$. 
Note that the two distinguished curves are then homotopic in $S$ relative to $\partial S$. 
The distinguished curves and the two components of their complement in the boundary of the quadrilateral are called the {\it edges} of the quadrilateral.

If $C$ is any arc in $S$ which is homotopic to one of the two distinguished edges of the quadrilateral, then we say that the quadrilateral has {\it core arc} $C$.

Any cylinder or quadrilateral in $S$ is equipped with a conformal structure induced from that of $S$. 
In this way, a cylinder in $S$ is conformally equivalent to a Euclidean annulus of the form $\{z\in\mathbb{C},\ r_1\leq \vert z\vert\leq r_2\}$. 
Although such a Euclidean annulus is not unique, 
the quotient ${r_1}/{r_2}$ is well-defined up to the choice of the boundary component of the cylinder that is sent to the interior component of the Euclidean annulus, 
and it is called the {\it modulus} of the Euclidean annulus and of the given cylinder in $S$. 
(Sometimes a factor of ${1}/{2\pi}$ is introduced in the definition of the modulus, but this is irrelevant to the present paper.)

Likewise, any quadrilateral  in $S$ is conformally equivalent to a Euclidean rectangle of the form $\{z\in\mathbb{C}, \ 0\leq \Re(z)\leq r_1,  0\leq \Im(z)\leq r_2\}$, 
with the distinguished edges sent to the vertical edges of this Euclidean rectangle.  
Again, although such a Euclidean rectangle is not unique, the quotient ${r_1}/{r_2}$ is uniquely defined, and it is called the {\it modulus} of the Euclidean rectangle 
and of the given quadrilateral in $S$.

The modulus of a quadrilateral or of an annulus $A$ in $S$ is denoted by $m(A)$.

\begin{definition}[Geometric] 
The extremal length $E(\gamma)$ of an element $\gamma$ in $\mathcal{C}(S)$ 
(respectively in $\mathcal{B}(S)$) is equal to \[\inf_A \frac{1}{m(A)},\] 
where the infimum is taken over all cylinders (respectively quadrilaterals) in $S$ with core curve (respectively core arc) in the homotopy class $\gamma$.
\end{definition}

We shall use the fact that the extremal length function has a continuous extension from the sets of homotopy classes of arcs and curves to the space of measured foliations. 
We refer to Kerckhoff \cite{K} for the case of surfaces without boundary, and we note the following, concerning the definition of a measured foliation on a surface with boundary, which make the results apply for such surfaces:

A measured foliation on a surface $S$ with boundary is a measured on the interior of $S$, which extends to the boundary of $S$, such that furthermore the following holds:
\begin{itemize}
\item The double of any measured foliation on $S$ is a measured foliation on $S^d$.
\item The restriction to $S$ of any measured foliation on $S^d$ is a measured foliation on $S$. 
Here, the surface $S$ is considered as embedded in $S^d$ in such a way that its boundary $\partial S$ is in minimal position with respect to the foliation on $S^d$ 
(cf. \cite{FLP} p.76 for the notion of minimal position, called there ``quasi-transverse" position).
\end{itemize}

Kerckhoff's formula for the Teichm\"uller metric in terms of extremal length is generalized  as follows for surfaces with boundary:

\begin{theorem}\label{th:Kerckhoff} 
If $X$ and $Y$ are two conformal structures on $S$, then the Teichm\"uller distance between them is equal to 
\[
d_{T}(X,Y)= \frac{1}{2}\log \sup_{\gamma\in\mathcal{C}\cup\mathcal{B}}\frac{E_Y(\gamma)}{E_X(\gamma)}.
\]
\end{theorem}
  
\begin{proof}  
In the case of closed surfaces, the proof of this formula is in \cite{K} p. 36 
(with the supremum being taken over $\gamma$ in $\mathcal{C}$, since the set $\mathcal{B}$ is empty). 
For a surface of finite type without boundary, the same proof applies because the curves which are parallel to the punctures (that is, the peripheral curves), are not part of the same $\mathcal{C}$, and therefore they are not involved in the proof. 
For general surfaces of finite topological type, the formula is obtained by taking the double, and we now outline the proof.

We have $d_{T}(X,Y)=d_{T}(X^d,Y^d)$, since Teichm\"uller distances are realized by Teichm\"uller maps, and, as we recalled above, 
a Teichm\"uller map between the surfaces without boundary $X^d$ and $Y^d$ induces a map between $X$ and $Y$, 
with the same dilatation, which will also be the Teichm\"uller map between the two doubles, since Teichm\"uller maps are unique. 
(Note that we are using the same notation $d_T$ for the Teichm\"uller distance on $\mathcal{T}(S)$ and of $\mathcal{T}(S^d)$.)
    
By Kerckhoff's formula for surfaces without boundary, we have
 
\[d_{T}(X^d,Y^d)= \frac{1}{2}\log \sup_{\gamma\in\mathcal{C}(S^d)}\frac{E_{Y^d}(\gamma)}{E_{X^d}(\gamma)}\]
and the supremum is realized by the extremal length of the vertical foliation of the initial quadratic differential $q$ associated to the Teichm\"uller map 
between the conformal structures $X$ and $Y$ (see also \cite{K} p. 36). 
Since the quadratic differential and its vertical foliation are invariant by the natural involution of $X^d$, 
and since the Teichm\"uller map restricts to a Teichm\"uller map between the two surfaces $X$ and $Y$, the supremum 
 \[ \sup_{\gamma\in\mathcal{C}(S)\cup \mathcal{B}(S)}\frac{E_{Y}(\gamma)}{E_{X}(\gamma)}\]
is realized by the restriction on $X$ of the vertical foliation of the quadratic differential $q$. 
This implies 
 \[d_T(X,Y)\leq  \frac{1}{2}\log \sup_{\gamma\in\mathcal{C}(S)\cup\mathcal{B}(S)}\frac{E_Y(\gamma)}{E_X(\gamma)}.\]
But the last inequality must be an equality, since if there were an element $\gamma$ in $\mathcal{C}(S)\cup \mathcal{B}(S)$ satisfying  
 \[d_T(X,Y)<  \frac{1}{2}\log  \frac{E_Y(\gamma)}{E_X(\gamma)},\]
then, taking the doubles of the Teichm\"uller maps between $X$ and $Y$ and the double $\gamma^d$ of $\gamma$, 
we would have an element $\gamma^d$ in $\mathcal{C}(S^d)\cup\mathcal{B}(S^d)$ satisfying 
   \[d_T(X^d,Y^d)<  \frac{1}{2}\log \frac{E_Y(\gamma^d)}{E_X(\gamma^d)},\]
which by Kerckhoff's formula for the Teichm\"uller metric of surfaces without boundary is excluded.
\end{proof}

We now deduce from  Theorem \ref{th:Kerckhoff} another formula for the Teichm\"uller metric on the Teichm\"uller space of a surface with boundary, 
which instead of using the extremal lengths of closed curves and of arcs, uses only extremal lengths of arcs. 
The idea is to approximate curves with arcs. 
We start with the following lemma which is extracted from \cite{LPST1}:
  
\begin{lemma}
\label{lemma:sym}
Let $\beta$ be a component of $\partial S$, let $\alpha$
be a measured foliation on $S$ and let $\alpha^d$ be its double. 
Then, there exists a sequence of symmetric simple closed curves on $S^d$ converging to $\alpha^d$ in the topology of  
$\mathcal{PMF}(S^d)$ such that each of these simple closed curves intersects essentially in exactly two points the image of $\beta$ in $S^d$, 
and intersects no other component of the image of $\partial S$ in $S^d$.
\end{lemma}

From this we deduce the following:
 
\begin{theorem}\label{th:Kerckhoff2} 
If $X$ and $Y$ are two conformal structures on $S$, then the Teichm\"uller distance between them is equal to 
\[
d_{T}(X,Y)= \frac{1}{2}\log \sup_{\mathcal{B}}\frac{E_Y(\gamma)}{E_X(\gamma)}.
\]
\end{theorem}

\begin{proof} 
We use the continuous extension of the extremal length function on the space $\mathcal{MF}(S)$ of measured foliations. 
Let $\alpha$ be a simple closed curve on $S$.
Consider its double $\alpha^{d}$ in $S^d$. 
By Lemma \ref{lemma:sym}, there exists a sequence of symmetric weighted simple closed curves $(\gamma^d_{n})$ such that
$\lim_{n\to\infty}\gamma^d_{n}=\alpha^d$ in $\mathcal{MF}(S)$.  
Let $X^d$ and $Y^d$ be the doubles of the  conformal structures $X$ and $Y$. We can assume
the sequence of symmetric simple closed curves $(\gamma^d_{n})$ on $S^d$ to be equipped, as elements of $\mathcal{MF}(S^d)$, 
 with the counting measure. Then, we have
$$
\Big{|}\frac{E_{Y^d}(\gamma^d_{n})}{E_{X^d}(\gamma^d_{n})}-
\frac{E_{Y^d}(\alpha^d)}{E_{X^d}(\alpha^d)}\Big{|}\to0.
$$

We conclude that for any simple closed curve $\alpha$ in $S$, there exists a sequence of connected arcs $(\gamma_{n})$ in $\mathcal{B}$ 
such that
$$
\Big{|}\frac{E_{Y}(\gamma_{n})}{E_{X}(\gamma_{n})}-
\frac{E_{Y}(\alpha)}{E_{X}(\alpha)}\Big{|}\to0.
$$
  
This gives
$$
\sup_{\alpha\in\mathcal{B}\cup\mathcal{C}}\frac{E_{Y}(\alpha)}{E_{X}(\alpha)}
=
\sup_{\gamma\in\mathcal{B}}\frac{E_{Y}(\gamma)}{E_{X}(\gamma)},
$$
which concludes the proof.
\end{proof}

Let us study a little further the function 
\[(X,Y)\mapsto \frac{1}{2}\log \sup_{\mathcal{C}}\frac{E_Y(\gamma)}{E_X(\gamma)}.\]

First note that, in the case where the surface $S$ is a pair of pants, this function does not define a metric, since it can take negative values. 
To see this, we take two conformal structures $X$ and $Y$ on $S$ such that for each of the structures, 
the three boundary components have the same extremal length 
(for instance, we can take two hyperbolic structures such that the three boundary components of each one have the same lengths) 
and such that the value of the extremal length for the conformal structure $Y$ is less than the corresponding value for the structure $X$. 
It is clear in this case that $\displaystyle  \frac{1}{2}\log \sup_{\mathcal{C}}\frac{E_Y(\gamma)}{E_X(\gamma)}<0$.
Actually, one has more generally the following result, whose proof follows.

\begin{proposition}\label{prop:ineq}
Let $X$ be a Riemann surface with boundary. 
Then there exists a Riemann surface structure $Y\not= X$ on the same topological surface such that 
$$\frac{1}{2}\log \sup_{\mathcal{C}}\frac{E_Y(\gamma)}{E_X(\gamma)}\leq0.$$
\end{proposition}

In particular, this shows that the function $\displaystyle \frac{1}{2}\log \sup_{\mathcal{C}}\frac{E_Y(\gamma)}{E_X(\gamma)}$ on $\mathcal{T}(S)$ does not separate points.

\begin{remark} In \cite{Parlier}, Parlier proved an analogue of Proposition \ref{prop:ineq} in the case
the case where extremal lengths are replaced by hyperbolic lengths.
\end{remark}

Before proving this proposition, let us recall a few facts on Nielsen extensions.

Let $X$ be a Riemann surface with boundary.
Consider the unique complete hyperbolic metric associated to the conformal structure on the interior of $X$ by the Uniformization Theorem. 
This hyperbolic metric is complete and there exists, for each boundary component $C_i$, $i=1, \cdots,b$, of $X$,
a unique simple closed geodesic $C'_i$ which is freely homotopic to $C_i$ in $S$.
By cutting $X$ along the geodesics $C'_i, i=1, \cdots,b$, 
we obtain a Riemann surface $X_0$ and $b$ \textsl{funnels} adjacent to $C'_i, i=1, \cdots, b$. 
Here a \textsl{funnel} is a Riemann surface of type $(0,0,2)$, that is, homeomorphic to an annulus.
The hyperbolic surface with totally geodesic boundary $X_0 \subset X$ is called the \textsl{Nielsen kernel} of $X$, 
and $X$ is the \textsl{Nielsen extension} of $X_0$ (see \cite{Bers}). 
The Nielsen kernel is a deformation retract of $X$ and can be viewed as embedded in $X$.

The proof of Proposition \ref{prop:ineq} uses the following:

\begin{lemma}\label{lem:mono}
Let $X_1$ be the Nielsen extension of $X_0$. 
Then for any $\gamma\in \mathcal{C}$, $E_{X_{1}}(\gamma)\le E_{X_{0}} (\gamma)$.
\end{lemma}

\begin{proof} 
This is a consequence of a monotonicity property for extremal length, which follows trivially from the definition.
\end{proof}

For any $X,Y$ in $\mathcal{T}(S)$, set
$$
T(X,Y)=\frac{1}{2}\log\sup_{\gamma \in \mathcal{C}} \frac{E_Y(\gamma)} {E_X (\gamma)}
$$
and
$$
\overline{T}(X,Y)=\frac{1}{2}\log\sup_{\gamma \in \mathcal{C}}
\frac{E_X(\gamma)} {E_Y (\gamma)}.
$$

From Lemma \ref{lem:mono}, we readily get 
$$
T(X,X_n)=\frac{1}{2}\log\sup_{\gamma \in \mathcal{C}} \frac{E_{X_n}
 (\gamma)} {E_X (\gamma)}\le 0.
$$

For the Riemann surface $X$, let $X_1$ be the Nielsen extension of $X$ and let $X_{k+1}$ be the Nielsen extension of $X_k$, $k=1,2,\cdots$. 
In view  of the  canonical embeddings $X_k \hookrightarrow X_{k+1}$ one can define the Riemann surface $X_{\infty}=X_1 \cup X_2 \cup \cdots$ 
which is called the \textsl{infinite Nielsen extension} of $X$.

Bers \cite{Bers} proved the following.

\begin{lemma}[Bers]
\label{le:infinite_Nielsen} 
The infinite Nielsen extension $X_{\infty}$ is a Riemann surface of type $(g,p+b,0)$.
\end{lemma}

Let us recall the following result of Maskit \cite{Maskit} (which is easily extended from the case without boundary to the case with boundary by taking doubles).

\begin{proposition}[Maskit]
\label{prop:Maskit}
Let $X$ be a Riemann surface of finite topological type and let $\gamma\in\mathcal{C}$. 
We have
$$
\frac{l_X (\gamma)}{\pi}\le E_X(\gamma) \le \frac{l_X (\gamma)}
{2} e^{\frac{l_X (\gamma)} {2}}.
$$
\end{proposition}

From Lemma \ref{le:infinite_Nielsen}, we have $\lim_{n\to\infty}l_{X_n}(C'_i)=0, i=1,\cdots,b$. 
Combining with Lemma \ref{prop:Maskit}, we obtain
$$
\lim_{n\to\infty}\overline{T}(X,X_n)=\lim_{n\to\infty}\frac{1}{2}\log\sup_{\gamma \in \mathcal{C}} \frac{E_X
 (\gamma)} {E_{X_n} (\gamma)}=\infty.
$$

We say that two functions $f$ and $g$ defined over $\mathcal{T}(S)$ are called $(\lambda,C)$-\textsl{quasi-isometric} 
if there exist constants $C\geq 0$ and $\lambda\geq 1$ such that, for any $X,Y\in \mathcal{T}(S)$,
$$\frac{1}{\lambda}g(X,Y) - C\leq f(X,Y)\leq \lambda g(X,Y) + C.$$

From what precedes we get the following

\begin{theorem} 
The functions $T$ and $\overline{T}$ are not quasi-isometric on $\mathcal{T}(S)$.
\end{theorem}

\begin{question} 
$\ $
\begin{itemize}
\item Let $T'(X,Y)=\max \left(T(X,Y), \overline{T}(X,Y)\right)$, for $X, Y \in \mathcal{T}(S)$. 
Is $T'(X,Y)$  a metric on $\mathcal{T}(S)$ ?\\
\item  Does there exist $X,Y \in \mathcal{T}(S)$, such that $T(X,Y)<0$?
\end{itemize}
\end{question}

\section{The relative thick part of Teichm\"uller space and Mumford's compactness theorem}

The {\it extended mapping class group of $S$}, denoted by $\mathrm{MCG}(S)$, is the group of homotopy classes of homeomorphisms of this surface. 
This group acts naturally on the Teichm\"uller space $\mathcal{T}(S)$, and the quotient space by this action is the {\it moduli space} $\mathcal{M}(S)$ of $S$.  Elements of $\mathcal{M}(S)$ are therefore homotopy classes of hyperbolic surfaces without marking.

For any given $\epsilon >0$, we call the {\it $\epsilon$-thick part of moduli space} the subset of $\mathcal{M}(S)$ consisting of homotopy classes of unmarked hyperbolic surfaces satisfying the following two conditions:
\begin{enumerate}
\item The length of any simple closed geodesic is $\geq \epsilon$.
\item The length of any geodesic arc is $\geq \epsilon$.
\end{enumerate}
 
A theorem by D. Mumford \cite{Mumford} says that in the case of surfaces of finite topological type without boundary, 
the $\epsilon$-thick part of moduli space 
(which in this case is simply the subset of $\mathcal{M}(S)$ consisting of homotopy classes of surfaces satisfying Condition (1)) 
is compact. 
The proof of this theorem that is given in \cite{Hubbard} can be adapted without difficulty to the case of hyperbolic surfaces with boundary, 
if we include Condition (2) on geodesic arcs joining boundary components.
For future reference, we state this as follows:

\begin{theorem}[Mumford's compactness theorem for surfaces with boundary]
\label{th:Mumford} 
For any topologically finite type surface $S$ and for any $\epsilon >0$, the $\epsilon$-thick part of the moduli space of $S$ is compact.
\end{theorem}

 In the paper \cite{LPST1}, given two real numbers $\epsilon$ and $\varepsilon_0$ satisfying $0<\epsilon \leq \varepsilon_0$, we considered the subspace of the Teichm\"uller space $\mathcal{T}(S)$, which we called the {\it $\varepsilon_0$-relative $\epsilon$-thick part of Teichm\"ulller space}, consisting of equivalence classes of hyperbolic metrics satisfying the following two conditions:
\begin{enumerate}[(a)]
 \item The length of any element in $\mathcal{C}$ is $\geq \epsilon$.
 \item The length of any boundary component of $S$ is $\leq \varepsilon_{0}$.
\end{enumerate}
We have the following

\begin{proposition}
For any real numbers $\epsilon$ and $\varepsilon_0$ satisfying $0<\epsilon \leq \varepsilon_0$, 
the natural image of the $\varepsilon_0$-relative $\epsilon$-thick part of Teichm\"ulller space in the moduli space $\mathcal{M}(S)$ 
is contained in some $\epsilon'$-thick part of moduli space, for some $\epsilon'>0$ depending only on $\epsilon$ and $\varepsilon_0$.
\end{proposition}

\begin{proof}
We showed in \cite{LPST1} (Lemma 3.4) that for all $X$ in the $\varepsilon_0$-relative $\epsilon$-thick part of $\mathcal{T}(S)$, 
to any geodesic arc $\beta$  on $S$, we can associate a simple closed geodesic $\alpha$ on $S$ such that the ratio 
$l_X(\beta) / l_X(\alpha)$ is bounded above and below by positive constants that depends only on $\epsilon$ and $\varepsilon_0$. 
In particular, if the length of the simple closed geodesic $\alpha$ is bounded below by $\epsilon$, the length of the arc $\beta$ is also bounded below by a constant that depends only on $\epsilon$ and $\varepsilon_0$. 
This gives the desired result.
\end{proof}

\section{Length spectrum weak metrics}\label{section:weak}

We recall that a  {\it weak metric} on a set $M$ is a function
$\delta:M\times M\to [0,\infty)$   satisfying
\begin{enumerate}[(a)]
 \item for $x$ and $y$ in $M$, $\delta(x,y)=0\iff x=y$;
\item  for
all $x$, $y$ and $z$ in $M$, $\delta(x,y)+\delta(y,z)\geq \delta(x,z)$.
\end{enumerate}
and we say that the weak metric $\delta$ is {\it asymmetric} if furthermore
 \medskip
 
\noindent (c)  there exist $x$ and $y$ in $M$ satisfying 
$\delta(x,y)\not=\delta(y,x)$.
\medskip

In the paper \cite{LPST1}, we considered the following two functions on $\mathcal{T}(S)\times \mathcal{T}(S)$, where $X$ and $Y$ are represented by hyperbolic structures:

\begin{equation}
\label{eq:def1} 
d(X,Y) = \log \sup_{\alpha \in  \mathcal {C}\cup  \mathcal {B}}\frac{l_Y(\alpha)}{l_X(\alpha)},
\end{equation}

\begin{equation}
\label{eq:def2}
\overline{d}(X,Y) = \log \sup_{\alpha \in \mathcal {C}\cup  \mathcal {B}}\frac{l_X(\alpha)}{l_Y(\alpha)}.
\end{equation}

The  functions in (\ref{eq:def1}) and (\ref{eq:def2}) are analogues, for surfaces with boundary, of asymmetric weak metrics introduced by Thurston in \cite{Thurston1998} for surfaces of finite type without boundary. 
We showed that these functions are asymmetric weak metrics on $\mathcal{T}(S)$ (Proposition 2.9 of \cite{LPST1}, where the point is to prove the separation property). 
We also showed that these weak metrics  can be expressed using suprema over the set $\mathcal{B}$ only. 
More precisely, we proved the following (\cite{LPST1} Proposition 2.12):

For every $X$ and $Y$ in $\mathcal{T}(S)$, we have
\begin{equation}\label{eq:def3}
d(X,Y)=\log\sup_{\gamma\in\mathcal{B}}\frac{l_{Y}(\gamma)}{l_{X}(\gamma)}.
\end{equation}

We also proved the following equalities (Corollary 2.8 of \cite{LPST1})

\begin{equation}\label{eq:double}
d(X,Y)=d(X^d,Y^d),\quad \overline{d}(X,Y)=\overline{d}(X^d,Y^d).
\end{equation}

The equalities in (\ref{eq:double}) are useful for obtaining results on the weak metrics $d$ and $\overline{d}$  for a surface $S$ with boundary from results on the corresponding weak metrics on the double of $S$, which is a surface without boundary.

We also studied the following symmetrization of the weak metrics $d$ and $\overline{d}$:

\begin{equation}\label{eq:def5}
\delta_L(X,Y) = \log\max \left(\sup_{\gamma \in \mathcal{C}\cup\mathcal{B}}\frac{l_Y(\gamma)}{l_X(\gamma)},
\sup_{\gamma \in \mathcal {C}\cup \mathcal{B}}\frac{l_X(\gamma)}{l_Y(\gamma)}\right).
\end{equation}
 which, by (\ref{eq:def3}) and (\ref{eq:def4}), can also be expressed as
\begin{equation}\label{eq:def6}
 \delta_{L}(X,Y)=\log\max\left(
 \sup_{\gamma\in\mathcal{B}}\frac{l_{Y}(\gamma)}{l_{X}(\gamma)},
 \sup_{\gamma\in\mathcal{B}}\frac{l_{X}(\gamma)}{l_{Y}(\gamma)}
 \right)=\max\left(d(X,Y),\overline{d}(X,Y)\right).
\end{equation}
 
From (\ref{eq:double}), we  immediately deduce the following, for all $X,Y$ in $\mathcal{T}(S)$:
 
 \begin{equation}\label{eq:def7}
 \delta_L (X,Y)=\delta_L (X^d, Y^d).
 \end{equation}

In the same paper \cite{LPST1}, we also considered the following symmetric function on the product $\mathcal {T}(S)\times \mathcal {T}(S)$:

\begin{equation}\label{eq:def8}
d_L(X,Y) = \log\max \left(\sup_{\gamma \in \mathcal {C}}\frac{l_X(\gamma)}{l_Y(\gamma)},
\sup_{\gamma \in \mathcal
{C}}\frac{l_Y(\alpha)}{l_X(\alpha)}\right).
\end{equation}

For surfaces of finite type without boundary, the functions $\delta_L$ and $d_L$ obviously coincide, 
and they define the so-called {\it length spectrum metric} on Teichm\"uller space. 
This metric was originally considered by Sorvali \cite{Sorvali},  
and it has been studied by several authors, see e.g. \cite {Liu2001} and \cite{CR}.
 
For surfaces of finite topological type with nonempty boundary, $\delta_L$ and $d_L$ are distinct, but they both are metrics. 
In the paper \cite{LPST1}, we gave a comparison between these metrics in the thick part of Teichm\"uller space.

Note the formal analogy between Kerckhoff's Formula adapted to surfaces with boundary (Theorem \ref{th:Kerckhoff}) and Version (\ref{eq:def1}) of the length spectrum asymmetric metric.
But since the apparently asymmetric formula in Theorem \ref{th:Kerckhoff} gives a genuine metric, 
Kerckhoff's formula is rather to be compared with the symmetrized length spectrum metric (\ref{eq:def6}). 
Our goal in the next section is to give a precise comparison between these two metrics. 
  
\medskip
  
\section{Comparison of length spectrum metrics with the Teichm\"uller metric}\label{section:comparison}
 
In this section, we use the preceding results to derive for surfaces with boundary some results which are already known for surfaces without boundary.
 
The first result is an immediate consequence of Wolpert's inequality (Theorem \ref{th:Wolpert} above):

\begin{theorem} 
\label{th:Wolpertcorollary}
For any surface $S$ of topologically finite type and for any $X$ and $Y$ in $\mathcal{T}(S)$, we have $\delta_L(X,Y)\leq 2d_T(X,Y)$.
\end{theorem}
  
The second result follows from a result obtained in \cite{Papado-Th2007} for surfaces without boundary.
\begin{theorem}
\label{th:co} 
Let $S$ be a surface of topologically finite type. 
The  weak metrics $d$, $\overline{d}$ and $\delta_L$ on $\mathcal{T}(S)$ are complete.
\end{theorem}
 
The result for the weak metrics $d$ and $\overline{d}$  follows from the relations (\ref{eq:double}) and the completeness of the corresponding weak metrics for the Teichm\"uller spaces of surfaces of finite type and without boundary (see \cite{Papado-Th2007}). 
The result for the metric $\delta_L$ follows from the definition, $\delta_L=\max (d,\overline{d})$.\\

The following is an adaptation of the result obtained by Choi \& Rafi in \cite{CR} for surfaces with boundary.
\begin{theorem} \label{theorem:thick}
Let $S$ be a surface of topologically finite type. 
For any $\epsilon$  and $\varepsilon_0$ satisfying $\varepsilon_0\geq \epsilon>0$ , there is a constant $D$ depending on
$\epsilon$ and $\varepsilon_0$ such that for any $X$ and $Y$ in the $\varepsilon_0$-relative $\epsilon$-thick part of $\mathcal {T}(S)$, we have
\begin{equation}\label{eq:CR1}
d(X, Y) - D \leq d_{T}(X, Y) \leq d(X, Y) + D,
\end{equation}
\begin{equation}\label{eq:CR2}
\overline{d}(X, Y) - D \leq d_{T}(X, Y) \leq \overline{d}(X, Y) + D,
\end{equation}
\begin{equation}\label{eq:CR3}
\delta_L(X, Y) - D \leq d_{T}(X, Y) \leq \delta_L(X, Y) + D
\end{equation}
 \end{theorem}

\begin{proof}  
We first prove the right hand side of Inequality (\ref{eq:CR1}).
 
Consider the function 
\begin{equation}\label{eq:CD5}
\displaystyle (X,\gamma)\mapsto \frac{E_X(\gamma)}{l_X^2(\gamma)},
\end{equation}
defined for $X$ in the $\varepsilon_0$-relative $\epsilon$-thick part of $\mathcal{T}(S)$ and $\gamma$ in $\mathcal{PMF}(S)$.
From the definition of the length of a measured lamination (see \cite{K}), we have the homogeneity property
$E_X(t\gamma) =t^2E_X(\gamma)$ and $l_X^2(t\gamma) = t^2l_X^2(\gamma)$ for all $t>0$. 
Since $\mathcal{PMF}(S)$ is compact and since the image of the  $\varepsilon_0$-relative $\epsilon$-thick part of $\mathcal{T}(S)$ in $\mathcal {M}(S)$ is compact (Theorem \ref{th:Mumford}), the values taken by the function (\ref{eq:CD5}) are uniformly bounded from above, that is, there exists a positive constant $C$ which depends only on $\epsilon$ and $\varepsilon_0$ such that

\begin{equation}\label{eq:compact}
\frac{E_X(\gamma)}{l_X^2(\gamma)} \leq C.
\end{equation}

It follows from (\ref{eq:extremal}) and (\ref{eq:compact}) that for any $X, Y$ in the $\varepsilon_0$-relative $\epsilon$-thick part of $\mathcal {T}(S)$, we have
$$\frac{E_Y(\gamma)}{E_X(\gamma)} \leq \frac{2\pi
|\chi(S)|}{l_X^2(\gamma)}\, C\, l_Y^2(\gamma).$$ 
Therefore,
\begin{eqnarray*}
d_{T}(X, Y) = \frac{1}{2} \log \sup_{\gamma \in \mathcal {C}(S)}
{\frac{E_Y(\gamma)}{E_X(\gamma)}} & \leq &\frac{1}{2}\log
\sup_{\gamma \in \mathcal {C}(S)}\frac{l_Y^2(\gamma)}{l_X^2(\gamma)}
+ \frac{1}{2}  \log{(2\pi |\chi(S)|\, C)}\\
& = & \log \sup_{\gamma \in \mathcal
{C}(S)}\frac{l_Y(\gamma)}{l_X(\gamma)}
+  \frac{1}{2} \log{(2\pi |\chi(S)|\, C)}\\
& = & d(X, Y) +  \frac{1}{2} \log{(2\pi |\chi(S)|\, C)}.\\
\end{eqnarray*}

From this, we obtain the right hand side inequality in (\ref{eq:CR1}).

The proof of the right hand side inequality in (\ref{eq:CR2}) is similar.

We can use (\ref{eq:extremal}) and (\ref{eq:compact}) again to show that
$$\frac{E_Y(\gamma)}{E_X(\gamma)} \geq \frac{l_X^2(\gamma)}{2\pi
|\chi(S)|Cl_Y^2(\gamma)},$$ which gives
$$d_{T}(X, Y) \geq 2\overline{d}(X, Y) - D\geq \overline{d}(X, Y) - D$$
for some $D$ depending on $\varepsilon_0$ and $\epsilon$.
Since $d_T$ is symmetric, we also have 
$$d_{T}(X, Y) \geq  d(X, Y) - D.$$

Inequalities (\ref{eq:CR3}) follow from the two previous ones and from the definition $\delta_L(X, Y) =
\max \left(d(X, Y), \overline{d}(X, Y)\right)$. 
\end{proof}

Threorem \ref{theorem:thick} gives in particular a comparison between the functions $d$ and $d_T$ for pairs of points in an $\varepsilon_0$-relative $\epsilon$-thick part of Teichm\"uller space. The following theorem gives information on pairs of points that are not necessarily in an $\varepsilon_0$-relative $\epsilon$-thick part of the space.

\begin{theorem}\label{th:LL}  
Let $S$ be a surface of topologically finite type and $(X_n)_{n\in\mathbb{N}}$ be a sequence of elements in $\mathcal{T}(S)$. 
Then,
  \[\lim_{n \to \infty}d_{T}(X_n, X_0)= \infty\iff \lim_{n \to \infty}d(X_n, X_0)= \infty.\]
 \end{theorem}

\begin{proof} 
Theorem \ref{th:Wolpertcorollary} implies that if $\lim_{n \to\infty}d(X_n, X_0)=\infty$, then we also have 
$\lim_{n \to\infty}d_{T}(X_n, X_0)=\infty$.\\

For the converse, suppose that $\lim_{n \to  \infty}d_{T}(X_n, X_0)=\infty$. 
First, consider the following two special cases:
\medskip

(1) The sequence $(X_n)_{n\in\mathbb{N}}$ stays in some $\varepsilon_0$-relative $\epsilon$-thick of
$\mathcal{T}(S)$ for some $\varepsilon_0\geq \epsilon>0$. 
In this case, by Theorem \ref{theorem:thick} (\ref{eq:CR1}), $\lim_{n \to\infty}d(X_n, X_0)=\infty$.
\medskip

(2) The sequence $(X_n)_{n\in\mathbb{N}}$ leaves any $\varepsilon_0$-relative $\epsilon$-thick part of $\mathcal{T}(S)$. 
In this case it follows from the definition of $d$ that $\lim_{n \to\infty}d(X_n, X_0)=\infty$.
\medskip
Now we discuss the  general case.  
We reason by contradiction.
Suppose that the sequence $d(X_n, X_0)$ does not go to infinity as $n$ goes to infinity. 
Then there exists a subsequence $(X_{n_k})_{k\in\mathbb{N}}$ of $(X_n)$ such that the set $\{d(X_{n_k},X_0), k\in\mathbb{N}\}$ is bounded. 
From the second special case above, if the subsequence $(X_{n_k})$ 
leaves any $\varepsilon_0$-relative $\epsilon$-thick part of $\mathcal{T}(S)$, we get a contradiction.
Therefore, up to taking a subsequence again, we can assume that there exists an $\varepsilon_0$-relative $\epsilon$-thick part of 
$\mathcal{T}(S)$ which contains the sequence of $(X_{n_k})$.
But the first case above shows that this yields a contradiction.\\
\end{proof}

\begin{theorem}
\label{th:top}
Let $S$ be a surface of topologically finite type and let $(X_n)_{n\ge 0}$ be a sequence of elements in $\mathcal{T}(S)$.
Then,
$$
\lim_{n\to\infty}d_T(X_n, X_0)=0\ \textrm{if and only if}\
\lim_{n\to\infty}\delta_L(X_n, X_0)=0;
$$
\end{theorem}

\begin{proof} 
We know that, for any $X,Y \in \mathcal{T}(S)$, 
$d_T(X, Y)=d_T(X^d,Y^d)$ and $\delta_L(X, Y)=\delta_L(X^d, Y^d)$. 
By the topological equivalence proved in \cite{Liu1999}, \cite{Papado-Th2007},
we have
$$
\lim_{n\to\infty}d_T(X^d_n, X^d_0)=0\ \textrm{if and only if}\ 
\lim_{n\to\infty}\delta_L(X^d_n, X^d_0)=0.
$$
This proves the theorem.
\end{proof}

We proved in \cite{LPST1} that for any surface $S$ of topologically finite type,  
  $d$, $\overline{d}$, $d_L$ and $\delta_L$ induce the same topology on $\mathcal {T}(S)$. Putting this result with Theorem \ref{th:top} we get the following

\begin{theorem}
Let $S$ be a surface of topologically finite type. 
The weak metrics $d$, $\overline{d}$, and the metrics $d_L$, $\delta_L$ and $d_T$ induce the same topology
on $\mathcal {T}(S)$.
\end{theorem}

\end{document}